\documentclass[12pt]{article} % Revision 8
\usepackage{amsfonts}
\usepackage{amsmath}

\title{A quadratic lower bound for subset sums}
\author{
   Matt DeVos
 \and
   Luis Goddyn\thanks{Supported by a Canada NSERC Discovery Grant}
 \and
   Bojan Mohar\thanks{Supported in part by the Slovenian ARRS Research
     Grant P1--0297 and in part by an NSERC Discovery Grant and
     CRC program.}
     \thanks{On leave from IMFM \& FMF, Department of
     Mathematics, University of Ljubljana, Ljubljana, Slovenia.}
 \and
   Robert \v{S}\'{a}mal\thanks{Supported by PIMS postdoctoral fellowship.}
   \thanks{On leave from Institute for Theoretical
     Computer Science (ITI), Charles University, Prague, Czech Republic.}\\
  \\
  {Department of Mathematics}\\
  {Simon Fraser University}\\
  {Burnaby, B.C. V5A 1S6} \\
  email: {\tt \{mdevos,goddyn,mohar,rsamal\}@sfu.ca}
}

\date{}
\newtheorem{theorem}{Theorem}[section]
\newtheorem{lemma}[theorem]{Lemma}
\newtheorem{corollary}[theorem]{Corollary}

\newtheorem{observation}[theorem]{Observation}

\newenvironment{proof}{\par\medskip\noindent{\bf Proof: }}
  {\hskip 2cm\unskip\hbox{}\hfill$\Box$\par\bigskip}
\newenvironment{proofof}{\par\medskip\proofof}
  {\unskip\hfill$\Box$\par\bigskip}
\def\proofof #1{\noindent{\bf Proof of #1:\hskip 0.5em}}

\newcommand\stab{\mathop{\mathit{stab}}}
\newcommand\df{\mathop{\mathrm{def}}\nolimits}
\newcommand\Cay{\mathop{\mathrm{Cayley}}}
\def\en{\mathbb{N}}
\def\zet{\mathbb{Z}}
\def\ceil#1{\lceil #1 \rceil}
\def\floor#1{\lfloor #1 \rfloor}
\def\bigfloor#1{\left\lfloor #1 \right\rfloor}
\def\?#1{{\bf ??? #1 ???}}
\begin{document}
\baselineskip=17pt

\maketitle

\begin{abstract}
Let $A$ be a finite nonempty subset of an additive abelian group
$G$, and let $\Sigma(A)$ denote the set of all group elements
representable as a sum of some subset of $A$. We prove that
$|\Sigma(A)| \ge |H| + \frac{1}{64}|A \setminus H|^2$ where $H$ is
the stabilizer of $\Sigma(A)$.  Our result implies that $\Sigma(A) =
\zet / n \zet$ for every set $A$ of units of $\zet / n \zet$ with
$|A| \ge 8 \sqrt{n}$.  This consequence was first proved by
Erd\H{o}s and Heilbronn for $n$ prime, and by Vu (with a weaker
constant) for general $n$.
\end{abstract}

\leftline{{\bfseries Keywords:\enspace}
  subset sums, additive bases, abelian groups
}

\leftline{{\bfseries MSC:\enspace}
11B13 % Additive bases
}

%%%%%%%%%%%%%%%%%%%%%%%%%%%%%%%%%%%%%%%
\section{Introduction}

%Throughout this paper we fix an additive abelian group $G$.
All groups considered in this paper are abelian, and we shall use
additive notation.  Let $G$~be such a group. If $A,B \subseteq G$,
then we let $A + B = \{ a + b : a \in A, b \in B \}$. If $g \in G$,
we let $g + A = A + g = \{g\} + A$, and we call any such set a
\emph{shift} of $A$.  The \emph{stabilizer} of $A$ is $\stab(A) = \{
g \in G : g + A = A \}$; note that this is a subgroup of~$G$.  We
define $\Sigma(A) = \{ \sum_{a \in A'} a : A' \subseteq A \}$, so
$\Sigma(A)$ is the set of group elements which can be represented as
sums of subsets of $A$. For any positive integer $n$, we let $\zet_n
= \zet/n\zet$.

In a lovely paper \cite{EH} which contains many of the ideas needed
in our proof, Erd\H{o}s and Heilbronn proved that $\Sigma(A) = G$
whenever $G \cong \zet_p$ for a prime $p$ and $A \subseteq G
\setminus \{0\}$ satisfies $|A| \ge 3 \sqrt{6p}$. They conjectured
that assuming $|A|\ge 2\sqrt p$ is sufficient; this was
confirmed by Olson~\cite{Ol} and further sharpened by
Dias da Silva and Hamidoune~\cite{SH}.

\begin{theorem} \label{olson_thm}
Let $p$ be a prime and let $A \subseteq \zet_p \setminus \{0\}$.
If\/ $|A| \le \floor{\sqrt{4p-7}}$, then $\Sigma(A) = \zet_p$.
\end{theorem}

To see that this theorem is essentially best possible, let $A
\subseteq \zet_p$ be the set $\{- \floor{\sqrt{p}},
\ldots, -1,1,\ldots,\floor{\sqrt{p}} \}$ and note that
$\floor{\frac{p}{2}} \not\in \Sigma(A)$. Such a strong conclusion does not
hold in general abelian groups, due to the existence of proper
nontrivial subgroups. For instance, if $H < G$ has $[G:H] = 3$ and
we take $A = H$, then $\Sigma(A) = H$ even though $A$ contains
one third of the elements in~$G$.  In cyclic groups, Vu found a
suitable assumption on $A$ which permits a similar conclusion.

\begin{theorem}[Vu \cite{Vu}] \label{vu_thm}
There exists a fixed constant $c$ so that $\Sigma(A) = \zet_n$ whenever
$A \subseteq \zet_n$ has size at least~$c \sqrt{n}$ and has the added
property that every number in $A$ is relatively prime with $n$.
\end{theorem}

The constant in this theorem is quite large.  It is derived from a
very deep theorem of Szemer\'{e}di and Vu~\cite{SV} on arithmetic
progressions in sumsets.  Our main theorem, which is quite
elementary by comparison, can be used to obtain Theorem \ref{vu_thm}
with a constant of $c = 8$.

Our main result gives a lower bound on $|\Sigma(A)|$, but before
introducing it, we shall pause to introduce Kneser's addition
theorem, an essential tool in our proof. Moreover, a simple
corollary of it gives a natural lower bound on $|\Sigma(A)|$ which
is of interest.

\begin{theorem}[Kneser \cite{Kn}]  \label{kneser}
Let $A_1,\ldots,A_m$ be finite nonempty subsets of~$G$.  If $H =
\stab(\sum_{i=1}^m A_i)$, then
$$ 
  \biggl| \sum_{i=1}^m A_i \biggr| \ge |H|(1-m) + \sum_{i=1}^m |A_i + H|.
$$
\end{theorem}

\begin{corollary}
Let $A \subseteq G$ and set $H = \stab(\Sigma(A))$.  Then
$$
|\Sigma(A)| \ge |H| + |H| \cdot |A \setminus H|.
$$
\end{corollary}

\begin{proof}
Let $A = \{a_1,\ldots,a_m\}$.  Then $\Sigma(A)
= \sum_{i=1}^m \{0,a_i\}$, and we obtain the desired bound by
applying Kneser's theorem to the right hand side of this equation.
\end{proof}

\bigskip

Our main theorem gives an alternative bound on $|\Sigma(A)|$ which
improves upon that from the previous corollary in the case when
$|H|$ is small.

\begin{theorem} \label{set_thm}
Let $A \subseteq G$ and set $H = \stab(\Sigma(A))$. Then
$$
    |\Sigma(A)| \ge |H| + \tfrac{1}{64}|A \setminus H|^2.
$$
\end{theorem}

As mentioned earlier, direct application of this result yields
Theorem~\ref{olson_thm} with a weaker constant and Theorem
\ref{vu_thm} with the stronger constant $c=8$.  To see this latter
implication, let $A \subseteq \zet_n$ have size $\ge 8 \sqrt{n}$,
assume it has the property that every element in $A$ is relatively
prime to $n$.  Suppose (for a contradiction) that $\Sigma(A) \neq
\zet_n$.  Then $H = \stab(\Sigma(A))$ is a proper subgroup of
$\zet_n$, so $A \cap H = \emptyset$ since every element in $A$
generates the entire group. But then our bound yields $|\Sigma(A)|
\ge |H| + \frac{1}{64} |A \setminus H|^2 > n$ --- a contradiction.

With some extra work we can improve our constant $1/64$ somewhat.
Indeed, it follows from our arguments that the same result holds
with a constant of ``almost'' $1/48$.  As far as we know,
Theorem~\ref{set_thm} may almost hold with $1/4$ in place of $1/64$:
it seems likely that $|\Sigma(A)| \ge \tfrac 14 |A \setminus H|^2 -
O(|A|)$. The extreme example we know of is essentially the same as
that mentioned earlier in connection with Olson's theorem. Namely,
if $A = \{-n,-(n-1),\ldots,n-1,n\} \subseteq \zet$. Then $|A| =
2n+1$, $H = \stab(\Sigma(A)) = \{0\}$ and $\Sigma(A) = \{ -
n(n-1)/2,\ldots,n(n-1)/2\}$ has size $n(n-1)+1$.

\bigskip

Theorem~\ref{set_thm} may be bootstrapped to give a bound on subsequence
sums. If ${\bf a}$ is a sequence of elements in $G$, we let $\Sigma({\bf a})$
denote the set of all sums of subsequences of ${\bf a}$. Note that if
${\bf a} = (a_1, \dots, a_n)$ and all the $a_i$'s are distinct then
$\Sigma({\bf a}) = \Sigma(\{a_1, \dots, a_n\})$; so subsequence sums
generalize the notion of subset sums.

%\footnote{that means elements of~$G/H \setminus \{H\}$}
If $H \le G$, we call any element of~$G/H \setminus \{H\}$
a \emph{nontrivial} $H$-coset of~$G$.
We let $\rho^j_H({\bf a})$ (for each $j \in \en$)
denote the number of nontrivial $H$-cosets of~$G$ 
which contain at least $j$ terms of ${\bf a}$.

\begin{theorem} \label{sequence_thm}
Let ${\bf a} = (a_1,\ldots,a_n)$ be a sequence
of elements in $G$, and let $H = \stab(\Sigma({\bf a}))$.
Then
$$
    |\Sigma({\bf a})|
       \ge |H| + \tfrac{1}{64} |H| \cdot
         \sum_{j \in \en} \bigl(\rho^j_H({\bf a})\bigr)^2.
$$
\end{theorem}

%%%%%%%%%%%%%%%%%%%%%%%%%%%%%%%%%%%%%%%
\section{Proofs}

The goal of this section is to prove our main results, Theorem
\ref{set_thm} and Theorem \ref{sequence_thm}.  In fact, these
theorems are easily seen to be equivalent, and our approach will be
to first prove Theorem \ref{set_thm} in the special case when $H =
\{0\}$, and then use this to prove the two main results in general.

Before we immerse ourselves into the details of the proof, let us
sketch our strategy. As in~\cite{EH}, the key goal is to show that
in every set $A \subseteq G$ with $|A| = 2(u+1)$ we can find a
subset~$B$ of size~$u+1$ such that $\Sigma(B)$ is large, provided
$\Sigma(A)$ has trivial stabilizer (Lemma~\ref{main_lem}).  To
establish this, we first use an inductive hypothesis to find a set
$B$ of size~$u$. Then we will try to find an element $c \in C = A
\setminus B$ such that by appending $c$ to~$B$, the size of
$S=\Sigma(B)$ grows significantly (thus maintaining our quadratic
bound). In other words, we want $\Delta_S(c) := |(S+c)\setminus S|$
to be large. Special cases of this task are dealt with in
Lemma~\ref{easy_bound} (if ``$S$ is small'') and~\ref{hard_bound}
(if ``$S$ is big''). In the work-horse of our proof,
Lemma~\ref{another_lemma}, we use these two to handle all possible
cases.

We also need to introduce a couple of definitions. If $G$~is a group
and $B \subseteq G$ then a \emph{(directed) Cayley graph~$\Cay(G,B)$} is
a graph with vertex-set~$G$ and with an arc $(g,g+b)$ for every $g \in G$
and $b \in B$.
If $B \subseteq G$ then we use $\langle B \rangle$ to denote the
subgroup of~$G$ generated by~$B$.

During the course of our proof we will often use Kneser's theorem
(Theorem~\ref{kneser}) and the following easy observations.

\begin{observation} \label{stab_obs}
We have $\stab(S) \le \stab(S+T)$ whenever $S,T \subseteq G$.

In particular, if $B \subseteq A$, then
$\stab(\Sigma(B)) \le \stab(\Sigma(A))$.
%because $\Sigma(Y) = \Sigma(X) + \Sigma(Y \setminus X)$.
\end{observation}

\begin{observation} \label{easy_obs}
If $A,B \subseteq G$ and $|A| + |B| > |G|$, then $A + B = G$.
\end{observation}

For every $S \subseteq G$ and every $x \in G$, we define
$\Gamma_S(x) = |(S+x) \cap S|$ and $\Delta_S(x) = |(S+x) \setminus
S|$.  Note that $\Gamma_S(x) + \Delta_S(x) = |S|$ and that
$\Delta_S(x) = \Delta_{G \setminus S}(x)$.  More interestingly, the
following observation shows that $\Delta_S$ is subadditive.

\begin{observation}[Erd\H{o}s, Heilbronn~\cite{EH}] \label{eh_obs}
If $x,y \in G$ then $\Delta_S(x+y) \le \Delta_S(x)
+\Delta_S(y)$.
\end{observation}

\begin{proof}
This is an immediate consequence of the following computation.
\begin{eqnarray*}
\Delta_S(x+y)
    & = &   |(S+x+y) \setminus S|                           \\
    &\le&   |(S+x+y) \setminus (S+y)| + |(S+y) \setminus S| \\
    & = &   |(S+x) \setminus S| + |(S+y) \setminus S|       \\
    & = &   \Delta_S(x) + \Delta_S(y)
\end{eqnarray*}
\end{proof}

If $Q,S \subseteq G$, we define the \emph{deficiency of $Q$
with respect to $S$} to be $\df_S(Q) = \min \{ |Q \cap S|, |Q\setminus S|\}$.

\begin{lemma} \label{easy_bound}
Let $C$, $S$ be finite subsets of a group $H$ such that
$\df_S(H) \le \frac{1}{2}|C|$.
Then $\frac{1}{|C|} \sum_{c \in C} \Delta_S(c) \ge \frac{1}{2} \df_S(H)$.
In particular, there exists $c \in C$ with
$\Delta_S(c) \ge \frac{1}{2} \df_S(H)$.
\end{lemma}

\begin{proof}
Recall that $\Delta_S(c) = \Delta_{H\setminus S}(c)$ for every~$c$.
Hence, after possibly replacing $S$ with $H \setminus S$,
we may assume that $\df_S(H) = |S|$.
Our lemma now follows from the inequalities below.
\begin{eqnarray*}
\sum_{c \in C} \Delta_S(c)
    & = &   |C| \cdot |S| - \sum_{c \in C} \Gamma_S(c)  \\
    &\ge&   |C| \cdot |S| - \sum_{h \in H} \Gamma_S(h) \\
    & = &   |C| \cdot |S| - |S|^2  \\
    &\ge&   \tfrac{1}{2} \, |C| \cdot |S|   \\
    & = &   \tfrac{1}{2} \, |C| \cdot \df_S(H)
\end{eqnarray*}
\end{proof}

\begin{lemma} \label{hard_bound}
Let $C$, $S$ be finite subsets of a group $H$ such that
$\df_S(H) \ge \frac{1}{2}|C|$ and $\langle C \rangle = H$.
Then there exists $c \in C$ with $\Delta_S(c) \ge \frac{1}{8}|C|$.
\end{lemma}

\begin{proof}
By possibly replacing $S$ with $H \setminus S$ we may assume that
$\df_S(H) = |S|$ and therefore $\tfrac 12 |C| \le |S| \le \tfrac{1}{2}|H|$.  
Now set $r = \floor{\frac{4|S|}{|C|}}$,
let $C^* = C \cup \{0\}$,
and let $D = \sum_{i=1}^r C^*$.  Put $K = \stab(D)$ and let
$t = |C^* + K|/|K|$, i.e., $t$ is the number of $K$-cosets
in~$H$ intersecting~$C^*$.
If $t \ge 2$, then by Kneser's addition theorem, we have
\begin{eqnarray*}
|D| &\ge&   r |C^* + K| - (r-1) |K| \\
    & = &   r (1 - \tfrac{1}{t})|C^* + K| + |K|    \\
    &\ge&   \bigl(\tfrac{4|S|-|C|}{|C|}\bigr) (1 - \tfrac{1}{t})|C| +
              \tfrac{1}{t} |C|    \\
    & = &   2|S| + \tfrac{t-2}{t}( 2|S| - |C| ) \\
    &\ge&   2|S|.
\end{eqnarray*}
If $t=1$, then $C \subseteq K$ and $C$ generates $H$, so we must
have $H = K = D$ and again we have $|D| \ge 2|S|$.  This brings us to
the following easy inequality:
$$
\sum_{d \in D} \Gamma_S(d)
    \le   \sum_{h \in H} \Gamma_S(h)
     =    |S|^2
    \le   \tfrac{1}{2} |D| \cdot |S|.
$$

It follows that there exists $d \in D$ with $\Gamma_S(d) \le \frac{1}{2}|S|$
and thus $\Delta_S(d) \ge \frac{1}{2}|S|$.  By construction, we may choose
elements $c_1,\ldots,c_n \in C$ with $n \le r$ so that
$d = \sum_{i=1}^n c_i$. Now, by the subadditivity of $\Delta_S$ we have
$\frac{1}{2}|S| \le \Delta_S(d) \le \sum_{i=1}^n \Delta_S(c_i)$ and
it follows that there exists an element $c \in C$ for which
$\Delta_S(c) \ge \frac{1}{2r}|S| \ge \frac{1}{8}|C|$ as desired.
\end{proof}

\begin{lemma} \label{another_lemma}
Let $A \subseteq G$ satisfy $|A| = 2u+2$ and ${\mathit
stab}(\Sigma(A)) = \{0\}$.  Let $\{B,C\}$ be a partition of $A$ with
$|B| = u$, put $S = \Sigma(B)$, and put $H = \langle C \rangle$.  If
$u \ge 16$ and $|H| \ge \frac{5}{256}u^2 + \frac{1}{4}u$, then one
of the following holds.
\begin{enumerate}
\item $|S| \ge \frac{1}{16}(u+1)^2$.
\item There exists $c \in C$ so that $\Delta_S(c) \ge \frac{1}{8}(u+1)$.
\end{enumerate}
\end{lemma}

\begin{proof} 
Define an $H$-coset $Q$ to be \emph{sparse} if $0 < |Q
\cap S| < \frac{1}{4}(u+1)$ and \emph{dense} if $|Q \setminus S| <
\frac{1}{4}(u+1)$.  If there is an $H$-coset $Q$ with $Q \cap S \neq
\emptyset$ which is neither sparse nor dense, then $\df_S(Q) \ge
\frac{1}{4}(u+1)$, so conclusion~2 follows by applying either
Lemma~\ref{easy_bound} or Lemma~\ref{hard_bound} to $C$ and an
appropriate shift of $Q \cap S$.  Thus, we may assume that every
$H$-coset which contains a point of $S$ is either sparse or dense.

If the sum of the deficiencies of the $H$-cosets (with respect to~$S$) is
at least $\frac{1}{4}(u+1)$, then by the averaging argument in
Lemma~\ref{easy_bound}, we find the existence of a $c \in C$ for which
$\Delta_S(c) \ge \frac{1}{8}(u+1)$ and conclusion~2 is
satisfied. Thus, we may assume that the sum of the deficiencies of
the $H$-cosets is at most $\frac{1}{4}(u+1)$. Since $|S| \ge u$ this
implies that there is at least one dense $H$-coset.

If $R$ is a dense $H$-coset, then it follows from
Observation~\ref{easy_obs} (and $|\Sigma(C)| \ge |C| \ge
\frac{1}{4}(u+1)$) that $\Sigma(C) + (R \cap S) = R$. Consequently,
if there are no sparse $H$-cosets, then $H \le {\mathit
stab}(\Sigma(C) + S) = \Sigma(A)$, which contradicts our
assumptions.  Thus, we may assume that there is at least one sparse
$H$-coset.  In particular, $S$ has nonempty intersection with at
least two $H$-cosets, so $S \not\subseteq H$.

If there exist four distinct dense $H$-cosets $Q_1$, \ldots, $Q_4$, then
we have the following:
\begin{eqnarray*}
|S|
    &\ge&   \sum_{i=1}^4 |S \cap Q_i|   \\
    & = &   4 |H| - \sum_{i=1}^4 \df_S(Q_i) \\
    &\ge&   \tfrac{5}{64}u^2 + u - \tfrac{1}{4}(u+1)  \\
    &\ge&   \tfrac{1}{16}(u+1)^2.
\end{eqnarray*}

Thus, we may assume that there are at most three dense $H$-cosets.
Now, for every $b \in B$, define
$S^+_b = b + \Sigma(B \setminus \{b\})$ and
$S^-_b = \Sigma(B \setminus \{b\})$.
Note that $S = S^+_b \cup S^-_b$.

\bigskip
\begin{narrower}
\noindent\textbf{Claim}
\begin{minipage}[t]{0.8\textwidth}
  If $R$ is a dense $H$-coset and $b \in B$,
  then either $R+b$ or $R-b$ is dense.
\end{minipage}

\begin{proof}
If $b \in H$, then the claim holds trivially,
so we may assume $b\not\in H$.  Let $d = \df_S(R)$ and
suppose (for a contradiction) that neither $R+b$ nor $R-b$ is dense.
Observe that $S \cap (R+b)$ contains $(S^-_b \cap R) + b$ and
$S \cap (R-b)$ contains $(S^+_b \cap R) - b$.  Suppose (without loss)
that $|S^-_b \cap R| \ge |S^+_b \cap R|$.  Then we have
\begin{eqnarray*}
\tfrac{1}{4}(u+1)
    & > &   \df_S(R) + \df_S(R+b)  \\
    &\ge&   d + |S^-_b \cap R|  \\
    &\ge&   d + \tfrac{1}{2}(|H|-d)  \\
    &\ge&   \tfrac{5}{512}u^2 + \tfrac{1}{8}u.
\end{eqnarray*}
This contradicts $u \ge 16$, thus establishing the claim.
\end{proof}
\end{narrower}

Let $W \subseteq G/H$ be the set of all dense $H$-cosets and set $w=
|W|$. We have that $1 \le w \le 3$ by our earlier arguments, but it
now follows from the claim (and $S \not\subseteq H$) that $w \ge 2$,
so $w \in \{2,3\}$.  For every $b \in G$, let $\Gamma_b$ be the
subgraph of $\Cay(G/H,b+H)$ induced by $W$. It follows from the
claim that $\Gamma_b$ has no isolated vertices whenever $b \in B
\setminus H$. Thus, every such $\Gamma_b$ is either a directed path
or a directed cycle.  If the graphs $\Gamma_b$ and $\Gamma_{b'}$
both have an edge with the same ends, then either $b'+H = b+H$ or
$b'+H = -b+H$.  It follows from this that either every $\Gamma_b$ is
a directed cycle, or every $\Gamma_b$ is a directed path; in the
latter case every pair of these paths have the same (unordered)
ends.
%all of these paths have the same two vertices as the ends
%ends of each of these paths are the same two vertices
%the same two vertices are the ends of each of these paths
%
If $\Gamma_b$ is a directed cycle for some $b \in B \setminus H$, then we
have $B \subseteq \langle H \cup \{b\} \rangle$ and we find that there are
no sparse $H$-cosets, contradicting our previous conclusions.

Thus, we may assume that every $\Gamma_b$ with $b \in B \setminus H$
is a directed path.  List the dense $H$-cosets $W_1,\ldots,W_w$ so
that every $\Gamma_b$ is a directed path with ends $W_1$ and $W_w$.
Setting $Q = W_2 - W_1$ we have that $W_1,\ldots,W_w$ is an
arithmetic progression in $G/H$ with difference $Q$.  Let $W_0 = W_1 - Q$
and $W_{w+1} = W_w + Q$; note that
$\{W_0,W_{w+1}\} \cap \{W_1,\ldots,W_w\} = \emptyset$.

Suppose first that $|B \setminus H| \ge w$. Choose $w$ distinct
elements $b_1,\ldots,b_w \in B \setminus H$ and for each of them
choose $\epsilon_i \in \{-,+\}$ so that $\epsilon_i b_i \in Q$. Now,
let $Z = W_1 \cap (\cap_{i=1}^w S^{-\epsilon_i}_{b_i})$ (in the
exponent we treat $\{+,-\}$ as a multiplicative group with identity
$+$).  It follows from our construction that $Z + \sum_{i=1}^w
\epsilon_i b_i \subseteq S \cap W_{w+1}$. For every $1 \le i \le w$
we have $(S^{\epsilon_i}_{b_i} \cap W_1) - \epsilon_i b_i \subseteq
W_0 \cap S$, so each $W_1 \cap S^{-\epsilon_i}_{b_i}$ contains all
but at most $|W_0 \cap S|$ points of $W_1 \cap S$. Thus, setting $d
= \df_S(W_1)$ we have the following inequalities (we use $|H| \ge
|C| > u+1$ and $|W_0 \cap S| < \frac{1}{4}(u+1)$).
\begin{eqnarray*}
\tfrac{1}{4}(u+1)
    & > &   \df_S(W_1) + \df_S(W_{w+1})   \\
    &\ge&   d + |Z| \\
    &\ge&   d + |H|- d - w |W_0 \cap S|  \\
    &\ge&   u+1 - \tfrac{w}{4}(u+1)
\end{eqnarray*}
However, this contradicts $w \in \{2,3\}$.  Thus $|B \setminus H| <
w$. But then, we must have $|B \setminus H| = w-1$, and we again
find that there are no sparse $H$-cosets, which contradicts our
previous conclusions.  This completes the proof. \end{proof}

Following Erd\H{o}s and Heilbronn~\cite{EH}, we define function
$L : \en \rightarrow \en$ by the following rule
$$
L(u) = \min_{ \begin{array}{c}
    \scriptstyle A \subseteq G \setminus \{0\} : |A| = 2u    \\
    \scriptstyle \stab(\Sigma(A)) = \{0\}
        \end{array} } \max_{B \subseteq A : |B| = u}  |\Sigma(B)|.
$$
We let $L(u) = \infty$ if no such set $A$ exists.  For every set~$B$
we have $\Sigma(B) \supseteq B \cup \{0\}$, so trivially $L(u) \ge u+1$.
Next we prove our main lemma which gives a better lower bound on~$L(u)$.

\begin{lemma} \label{main_lem}
$L(u) \ge \frac{1}{16}u^2$ for every $u \in \en$.
\end{lemma}

\begin{proof}
We proceed by induction on $u$.  Assume that the lemma holds for all
integers $\le u$ and let $A \subseteq G$ satisfy $|A| = 2(u+1)$ and
$\stab(\Sigma(A)) = \{0\}$. (If there is no such set, then we have
defined $L(u+1) = \infty$.)  We will show that there exists $B'
\subseteq A$ with $|B'| = u+1$ and $|\Sigma(B')| \ge
\frac{1}{16}(u+1)^2$.  It follows from our trivial bound $L(u)\ge
u+1$ that we may assume $u \ge 16$.

Apply the lemma inductively to obtain a set $B \subseteq A$ with
$|B| = u$ and $|\Sigma(B)| \ge \frac{1}{16}u^2$. Put $C = A
\setminus B$. To apply Lemma~\ref{another_lemma}, which is our aim,
we need a lower bound on the size of~$H = \langle C\rangle$. We do
this by estimating~$|\Sigma(C)|$. To this end, we apply the lemma
inductively twice more: choose a set $C_1 \subseteq C$ of size
$\ceil{\frac{u}{2}}$ with $|\Sigma(C_1)| \ge \frac{1}{64} u^2$ and
(since $2 \ceil{\frac{u}{4}} + \ceil{\frac{u}{2}}
  \le 2\frac{u+3}{4} + \frac{u+1}{2}
    = u + 2$)
a set $C_2 \subseteq C \setminus C_1$ of size $\ceil{\frac{u}{4}}$
with $|\Sigma(C_2)| \ge \frac{1}{256} u^2$. Put $C_3 = C \setminus
(C_1 \cup C_2)$. Now $\Sigma(C) = \Sigma(C_1) + \Sigma(C_2) +
\Sigma(C_3)$. Since $\Sigma(C)$ has trivial stabilizer 
(Observation~\ref{stab_obs}), Kneser's theorem gives the following
inequality (in the last inequality we use the trivial bound
$|\Sigma(X)| \ge |X| + 1$ if $0 \not \in X$).
\begin{eqnarray*}
|\Sigma(C)|
    & = & |\Sigma(C_1) + \Sigma(C_2) + \Sigma(C_3)| \\
    &\ge& |\Sigma(C_1)| + |\Sigma(C_2)| + |\Sigma(C_3)| - 2 \\
    &\ge& \tfrac{5}{256}u^2 + \tfrac{1}{4} u - 1.
\end{eqnarray*}
Let $S = \Sigma(B)$ and recall that $H = \langle C \rangle$.  Since
$\Sigma(C)$ has trivial stabilizer, $\Sigma(C) \subset H$, so $|H|
\ge \tfrac{5}{256}u^2 + \tfrac{1}{4}u$.  Since $u \ge 16$ by
assumption, we may apply Lemma~\ref{another_lemma} to deduce that
either $|S| \ge \frac{1}{16}(u+1)^2$, in which case we are finished,
or there exists $c \in C$ so that $\Delta_S(c) \ge
\frac{1}{8}(u+1)$.  In the latter case we have
\begin{eqnarray*}
|\Sigma(B \cup \{c\})|  & = &   |S + \{0,c\}|   \\
    & = &   |S| + \Delta_S(c)   \\
    &\ge&   \tfrac{1}{16}u^2 + \tfrac{1}{8}(u+1)  \\
    & > &   \tfrac{1}{16}(u+1)^2.
\end{eqnarray*}
This completes the proof.
\end{proof}

\begin{lemma} \label{triv_stab_lem}
If $A \subseteq G$ satisfies $\stab(\Sigma(A)) = \{0\}$, then
$$
  |\Sigma(A)| \ge 1 + \tfrac{1}{64}|A \setminus \{0\}|^2.
$$
\end{lemma}

\begin{proof}
We may assume that $0\not\in A$ since this has
no effect on our bound.  Set $|A| = u$.  The lemma holds trivially
if $u \le 8$, so we may assume $u>8$. By the previous lemma we may
choose a subset $B \subseteq A$ of size $\floor{\frac{u}{2}}$
such that $|\Sigma(B)| \ge L( \floor{\frac{u}{2}}) \ge \frac{(u-1)^2}{64}$.
Let $C = A \setminus B$. Then, by Kneser's theorem we have
\begin{eqnarray*}
|\Sigma(A)|
    & = &   |\Sigma(B) + \Sigma(C)|  \\
    &\ge&   |\Sigma(B)| + |C| - 1    \\
    &\ge&   \tfrac{1}{64}u^2 - \tfrac{1}{32}u + \tfrac{u}{2} - 1    \\
    &\ge&   1 + \tfrac{1}{64}u^2.
\end{eqnarray*}
\end{proof}

Note that by recursively applying Lemma \ref{main_lem} we can
improve the constant $1/64$ to ``almost'' $1/48$: 
In the above proof, we have used that $|\Sigma(B)| \ge L(|B|)$
(together with Lemma~\ref{main_lem}) to bound $|\Sigma(B)|$.
On the other hand, for $|\Sigma(C)|$ we have used a straightforward
bound $|\Sigma(C)| \ge |C|+1$. We could instead use the same procedure
recursively with $C$ in place of~$A$. This yields
$$
  |\Sigma(A)| \ge \frac {1}{16} \Bigl( 
    \bigfloor{\frac{u}{2}}^2 +
    \bigfloor{\frac{u}{4}}^2 +
    \bigfloor{\frac{u}{8}}^2 + \cdots \Bigr) 
  = \frac{1}{48} u^2 - O(u) \,.
$$

Next we prove our main theorem for sequences.

\bigskip
\begin{proofof}{Theorem~\ref{sequence_thm}}
For $Q \in G/H$ we let $c(Q) = |\{ 1 \le j \le n : a_j \in Q\}|$.
Let further $s$ be the maximum of $c(Q)$ for nontrivial cosets $Q$
(i.e., $Q \in G/H \setminus \{H\}$).
Finally, put $A_j = \{Q : c(Q) \ge j\}$ and note that
$|A_j| = \rho_H^j({\bf a})$.  By applying Kneser's theorem and
Lemma~\ref{triv_stab_lem} in the quotient group $G/H$, we have
\begin{eqnarray*}
|\Sigma( {\bf a} )|
    & = &   |H| \cdot \Bigl|\sum_{j=1}^{s} \Sigma(A_j)\Bigr|            \\
    &\ge&   |H| \cdot \Bigl(\sum_{j=1}^s |\Sigma(A_j)| - s + 1 \Bigr)   \\
    &\ge&   |H| + \tfrac{1}{64} |H| \cdot
      \sum_{j \in \en} \bigl(\rho^j_H({\bf a})\bigr)^2
\end{eqnarray*}
which completes the proof.
\end{proofof}

Finally, we prove our main theorem for sets.

\bigskip
\begin{proofof}{Theorem~\ref{set_thm}}
Let $A = \{a_1,a_2,\ldots,a_n\}$ and put ${\bf a} = (a_1, a_2, \ldots, a_n)$.
By applying Theorem~\ref{sequence_thm} we have
$$
|\Sigma(A)|
     =   |\Sigma( {\bf a} )|
    \ge  |H| + \tfrac{1}{64} |H| \cdot
      \sum_{j \in \en} \bigl(\rho^j_H({\bf a})\bigr)^2.
$$
Since $A$ is a set, $\rho^{j}_H( {\bf a} ) = 0$ for every $j > |H|$.
Further $\sum_{j=1}^{|H|} \rho^j_H( {\bf a}) = |A \setminus H|$.
It follows from this (and the Cauchy--Schwarz Inequality) that
$|H| \cdot \sum_{j \in \en} (\rho^j_H({\bf a}))^2 \ge |A \setminus H|^2$.
Combining this with the above inequality yields the desired bound.
\end{proofof}

%\bibliographystyle{rs-amsplain}
%\bibliography{subsum}

\end{document}